\documentclass[12pt]{article}
\usepackage{graphics,amsmath,amssymb}
\usepackage{latexsym}
\usepackage{epsfig}
\usepackage{hyperref}
\usepackage{fancybox}
\usepackage{xcolor}

\def\epsilon{\varepsilon}

\def\phi{\varphi}

\newtheorem{theorem}{Theorem}[section]
\newtheorem{lemma}[theorem]{Lemma}

\newtheorem{definition}[theorem]{Definition}

\newtheorem{proposition}[theorem]{Proposition}
\newtheorem{remark}[theorem]{Remark}

\def\N{{\mathbb N}}
\def\C{{\mathbb C}}
\def\R{{\mathbb R}}

\newenvironment{Proof}{\removelastskip\par\medskip
\noindent{\em Proof.} \rm}{\penalty-20\null$\square$\par\medbreak}

\title{\bf   Trace regularity for biharmonic 
evolution 
\\
equations 
with 
Caputo derivatives}

\author{Paola Loreti
\thanks{Dipartimento di Scienze di Base e Applicate per l'Ingegneria,
Sapienza Universit\`a di Roma,
Via Antonio Scarpa 16, 00161 Roma (Italy); e-mail: 
$<$paola.loreti@uniroma1.it$>$ }
\and Daniela Sforza
\thanks{Dipartimento di Scienze di Base e Applicate per l'Ingegneria, 
Sapienza Universit\`a di Roma,
Via Antonio Scarpa 16, 00161 Roma (Italy); e-mail: 
$<$daniela.sforza@uniroma1.it$>$ }
\thanks{To appear in: Fract. Calc. Appl. Anal. Vol. 25 (2022), at https://www.springer.com/journal/13540}}

\begin{document}

\maketitle

\begin{abstract}
Our goal is to establish a hidden regularity result for solutions of time fractional Petrovsky system. The order $\alpha$ of the Caputo fractional derivative belongs to the interval $(1,2)$. We achieve such result for a suitable class of weak solutions.

\end{abstract}



\bigskip
\noindent

\section{Introduction}

%
%
%
%
%

Our aim  is to study 
the fractional Petrovsky system 
\begin{equation}\label{eq:weakpI}
\begin{cases}
\partial_t^{\alpha}u+\Delta^2 u=0
\qquad \mbox{in}\ (0,T)\times\Omega,
\\
u=\Delta u=0    \hskip1.4cm \mbox{on}\ (0,T)\times\partial\Omega,
\end{cases}
\end{equation}
where  $\partial_t^{\alpha}u$ is the Caputo derivative with $\alpha\in(1,2)$. 

A first open question we intend to investigate is the existence of weak solutions for problem \eqref{eq:weakpI} in appropriate Sobolev spaces. We also analyse if trace regularity holds for weak solutions. For  non fractional hyperbolic partial differential equations trace regularity is also known as hidden regularity, a concept introduced by J.L. Lions in \cite{Lio2}, see \cite{LasTri} too.

In the non fractional setting the Petrovsky system
models vibrations of   beams and plates.
One can consider plate equations with different boundary conditions. Nevertheless, in the literature two classes of boundary conditions  are mainly discussed.
 The first one involves the function and its 
normal derivative on the boundary  (Dirichlet-Neumann)  and it describes the physical
model of  vibration of guided beams and plates. The second type of boundary conditions  includes the function and its Laplacian on the boundary: it depicts  the physical model of vibration of  hinged  beams and plates, see 
\cite {JL-Sia} and also \cite {K}.
Both problems have been studied in control theory to get the exact controllability by acting on (eventually on a part of) the boundary, see
\cite {Lio3}.

Nowadays, fractional calculus has attracted an increasing interest  in modeling evolution equations.  We refer to  Mainardi  and Gorenflo \cite {MG} for a survey devoted to the  theory of relaxation processes governed by linear differential equations of fractional order,   where the fractional derivatives are intended both in the sense of  Riemann-Liouville and in the sense of Caputo. 
To our knowledge,  the  time-fractional Petrovsky system has not
been  studied yet.  This paper would be a first attempt to solve some questions about trace regularity. 

To begin with, a difficulty
is an appropriate definition for  weak solutions. For a general analysis concerning the type of regularity required for solutions to fractional differential equations, see \cite{St}.

The definition of $H^2$-solutions is given in \cite{LS4}:
a $H^2$-solution of the fractional boundary value problem \eqref{eq:weakpI} is a
function $u$ belonging to
$C([0,T];H^2(\Omega)\cap H^1_0(\Omega))$ with $u_t\in L^2(0,T;L^2(\Omega))\cap C([0,T];D(\Delta^{-2\theta}))$,  for some $\theta\in(0,1)$, 
and satisfying, for any $v\in H^2(\Omega)\cap H^1_0(\Omega)$, 
\begin{equation*}
\int_\Omega I^{2-\alpha}\big(u_t(\cdot,x)-u_t(0,x)\big)v(x) \ dx\in C^1([0,T])
\end{equation*} 
and 
\begin{equation*}
\frac{d}{dt}\int_\Omega I^{2-\alpha}\big(u_t(t,x)-u_t(0,x)\big)v(x) \ dx+\int_\Omega\Delta u(t,x)\Delta v(x) \ dx=0
\qquad t\in (0,T)
\,.
\end{equation*}
The expression $H^2$-solution is suggested by the analysis of the stationary case given in  \cite {NS}. 

As a consequence of the existence and uniqueness result for $H^2$-solutions, the normal derivative of the solution on the boundary is well defined. 
Our purpose is to give a meaning to the normal derivative of weaker
solutions belonging to $H_0^1(\Omega)$. 

A function
$u\in C([0,T];H_0^1(\Omega))$ is called  a $H^1$-solution of \eqref{eq:weakpI} if, set $w=(-\Delta)^{-1} u$, we have
$w\in C([0,T];D((-\Delta)^{\frac32}))$, $w_t\in L^2(0,T;H^1(\Omega))\cap C([0,T];H^{-1}(\Omega))$,  $\partial_t^{\alpha} w\in L^2([0,T];L^2(\Omega))$
and 
\begin{equation*}
\int_\Omega \partial_t^{\alpha} w(t,x)v(x) \ dx-\int_\Omega\nabla\Delta w(t,x)\cdot\nabla v(x) \ dx=0 
\qquad v\in H^1_0(\Omega),\ t\in (0,T)
\,.
\end{equation*}
Thanks to this setting we are able to handle the question of 
well-posedness. Indeed, we establish the following result
\begin{theorem}\label{th:reg-l2}
Let $u_0\in H^1_0(\Omega)$ and $u_1\in H^{-1}(\Omega)$. Then the function
\begin{equation*}
u(t,x)=\sum_{n=1}^\infty\big[  \langle u_0,e_n\rangle E_{\alpha}(-\lambda_nt^\alpha) 
+\langle u_1,e_n\rangle_{H^{-1}} t E_{\alpha,2}(-\lambda_nt^\alpha)\big]e_n(x)
\end{equation*}
is the  $H^1$-solution  of \eqref{eq:weakpI}  satisfying the initial conditions
\begin{equation}\label{Aeq:cauchy10I}
u(0,\cdot)=u_{0},\quad
u_t(0,\cdot)=u_{1}.
\end{equation}
\end{theorem}
Next step is to prove the hidden regularity for $H^1$-solutions.
Although the result is similar to \cite[pag. 29] {K}, several technical steps
need a change in order to avoid integration by parts with respect
to the time.  

This difficulty also appears in \cite  {LSw} and we  borrow  from  \cite  {LSw} the way to deal with  some identities  useful in the proof of the hidden regularity result.

\begin{theorem}
For $u_0\in H^1_0(\Omega)$ and $u_1\in H^{-1}(\Omega)$, if $u$ is the $H^1$-solution  of \eqref{eq:weakpI}--\eqref{Aeq:cauchy10I}
then   for any $T>0$ we have
\begin{equation*}
\int_0^T\int_{\partial\Omega} \big|\partial_\nu u\big|^2d\sigma dt
\le C\Big(\| u_0\|^2_{H^1_0(\Omega)}+\|u_1\|^2_{H^{-1}(\Omega)}\Big)
\,,
\end{equation*}
for some constant $C=C(T)$ independent of the initial data.
\end{theorem}
The paper is organized as follows. 
In section 2 we give  a mathematical background 
on biharmonic operators and fractional derivatives.
In section 3 we introduce $H^1$ and $H^3$ solutions
for  fractional Petrovsky systems and we show existence theorems.
In section 4 we give a result of hidden regularity.

\section{Preliminaries}\label{s:pre}
Here we collect some notations, definitions and known results that we use to prove our main results.
\subsection{Biharmonic operators}\label{s:}
Let $\Omega\subset\R^N$,  $N\ge1$, be a bounded open set with sufficiently smooth  boundary denoted by $\partial\Omega$. As usual, we consider 
$
L^2(\Omega)
$
endowed with the inner product and norm defined by
\begin{equation*}
\langle u,v\rangle=\int_{\Omega}u(x)v(x)\ dx,
\qquad
\|u\|_{L^2(\Omega)}=\left(\int_{\Omega}|u(x)|^{2}\ dx\right)^{1/2}\qquad
u,v\in L^2(\Omega).
\end{equation*}
We define the operator $A$  in $L^2(\Omega)$ by
\begin{equation}\label{eq:defa}
\begin{split}
D(A)&= H^4(\Omega)\cap \{u\in H^3(\Omega): u=\Delta u=0\ \mbox{on}\ \partial\Omega\}
\\
(A u)(x)&=\Delta^2 u(x), 
\quad u\in D( A),
\  x\in\Omega.
\end{split}
\end{equation}
The spectrum of the operator $A$ consists of a sequence $\{\lambda_n\}$ tending to $+\infty$. Moreover, we assume that the eigenvalues $\lambda_n$ are all distinct numbers, and hence the eigenspace generated by $\lambda_n$  has dimension one. For special domains this assumption is fulfilled, e.g. if $\Omega$ is a ball in $\R^N$.

In addition we denote by $e_n$  the eigenfunctions of $A$  \big($A e_n=\lambda_n e_n$\big) that constitute an orthonormal basis of $L^2(\Omega)$.

The biharmonic operator given by \eqref{eq:defa} is 
self-adjoint, positive and 
 the domain $D(A)$ is dense in $L^2(\Omega)$.

The fractional powers $A^{\theta}$ are defined for $\theta>0$, see e.g. \cite{Pazy,Lunardi}. 
The domain $D(A^{\theta})$ of $A^{\theta}$ consists of those functions $u\in L^2(\Omega)$ such that
\begin{equation*}
\sum_{n=1}^\infty \lambda_n^{2\theta} |\langle u,e_n\rangle|^2<+\infty
\end{equation*}
and
\begin{equation*}
A^{\theta} u=\sum_{n=1}^\infty \lambda_n^{\theta} \langle u,e_n\rangle e_n,
\quad
u\in D(A^{\theta}).
\end{equation*}
Moreover $D(A^{\theta})$ is a Hilbert space with the norm given by
\begin{equation}\label{eq:norm-frac}
\|u\|_{D(A^{\theta})}=\|A^{\theta} u\|_{L^2(\Omega)}=\left(\sum_{n=1}^\infty \lambda_n^{2\theta} |\langle u,e_n\rangle|^2\right)^{1/2}
\quad
u\in D(A^{\theta}),
\end{equation}
and for any $0<\theta_1<\theta_2$ we have $D(A^{\theta_2})\subset D(A^{\theta_1})$.

In particular, $D(A^{\frac12})=H^2(\Omega)\cap H^1_0(\Omega)$ and 
$D(A^{\frac14})=H^1_0(\Omega)$ with the respective norms given by
\begin{equation}\label{eq:H^2_0}
\|u\|_{D(A^{\frac12})}=\left(\sum_{n=1}^\infty \lambda_n |\langle u,e_n\rangle|^2\right)^{1/2}=\|\Delta u\|_{L^2(\Omega)}
\qquad
u\in H^2(\Omega)\cap H^1_0(\Omega),
\end{equation}
\begin{equation}\label{eq:H^1_0}
\|u\|_{D(A^{\frac14})}=\left(\sum_{n=1}^\infty \lambda_n^{\frac12} |\langle u,e_n\rangle|^2\right)^{1/2}=\|u\|_{H^1_0(\Omega)}=\|\nabla u\|_{L^2(\Omega)}
\qquad
u\in H^1_0(\Omega).
\end{equation}
We also note that $D(A^{\frac34})=\{u\in H^3(\Omega): u=\Delta u=0\ \mbox{on}\ \partial\Omega\}$ and
\begin{equation}\label{eq:H^{-1}0}
\|u\|_{D(A^{\frac34})}=\left(\sum_{n=1}^\infty \lambda_n^{\frac32} |\langle u,e_n\rangle|^2\right)^{1/2}=\|\nabla\Delta u\|_{L^2(\Omega)}
\qquad
u\in D(A^{\frac34})
\,,
\end{equation}
see \cite[Lemma 1.7]{K}.
If we identify the dual $(L^{2}(\Omega))'$ with $L^{2}(\Omega)$ itself, then we have $D(A^{\theta})\subset L^{2}(\Omega)\subset(D(A^{\theta}))'$.
From now on we set 
\begin{equation}
D(A^{-\theta}):=(D(A^{\theta}))',
\end{equation}
whose elements are continuous linear functionals on $D(A^{\theta})$. If $u\in D(A^{-\theta})$ and $\varphi\in D(A^{\theta})$ the value $u(\varphi)$ is denoted  
by 
\begin{equation}\label{eq:duality}
\langle u,\varphi\rangle_{-\theta,\theta}:=u(\varphi)\,.
\end{equation}
According to the Riesz-Fr\'echet representation theorem $D(A^{-\theta})$ is a Hilbert space with the norm given by 
\begin{equation}\label{eq:norm-theta}
\|u\|_{D(A^{-\theta})}=\left(\sum_{n=1}^\infty \lambda_n^{-2\theta} |\langle u,e_n\rangle_{-\theta,\theta}|^2\right)^{1/2}
\qquad
u\in D(A^{-\theta})
\,.
\end{equation}
Moreover, 
\begin{equation*}
A^{-\theta} u=\sum_{n=1}^\infty \lambda_n^{-\theta} \langle u,e_n\rangle_{-\theta,\theta} e_n,
\quad
\|u\|_{D(A^{-\theta})}=\|A^{-\theta} u\|_{L^2(\Omega)}
\qquad
u\in D(A^{-\theta}),
\end{equation*}
and
for any $0<\theta_1<\theta_2$ we have $D(A^{-\theta_1})\subset D(A^{-\theta_2})$.
We also recall that 
\begin{equation}\label{eq:-theta}
\langle u,\varphi\rangle_{-\theta,\theta}=\langle u,\varphi\rangle
\qquad \mbox{for}\ u\in L^{2}(\Omega),\ \varphi\in D(A^{\theta}),
\end{equation}
see e.g. \cite[Chapitre V]{B}.
We note that $D(A^{-\frac14})=H^{-1}(\Omega)$ and
\begin{equation}\label{eq:H^{-1}}
\|u\|_{D(A^{-\frac14})}=\left(\sum_{n=1}^\infty \lambda_n^{-\frac12} |\langle u,e_n\rangle_{-\frac14,\frac14}|^2\right)^{1/2}=\|u\|_{H^{-1}(\Omega)}
\quad
u\in H^{-1}(\Omega)
\,.
\end{equation}

\subsection{Fractional derivatives}\label{s:}
\begin{definition}
For any $\beta>0$
we denote the Riemann--Liouville fractional integral operator of order $\beta$ by
\begin{equation}\label{eq:RLfi}
I^{\beta}(f)(t)=\frac1{\Gamma(\beta)}\int_0^t (t-\tau)^{\beta-1}f(\tau)\ d\tau, 
\qquad
f\in L^1(0,T),
\ \mbox{a.e.}\ t\in(0,T), 
\end{equation}
where $T>0$ and $\Gamma (\beta)=\int_0^\infty t^{\beta-1}e^{-t}\ dt$ is the Euler Gamma function. 
\end{definition}
The Caputo fractional derivative of order $\alpha\in(1,2)$ is given by
\begin{equation}\label{eq:der-frac}
\partial_t^{\alpha}f(t)
=\frac1{\Gamma(2-\alpha)}\int_0^t (t-\tau)^{1-\alpha}\frac{d^2f}{d\tau^2}(\tau)\ d\tau\,.
\end{equation}
By means of the Riemann--Liouville  integral operator $I^{2-\alpha}$ we can write
\begin{equation}\label{eq:der-fracI}
\partial_t^{\alpha}f(t)
=
I^{2-\alpha}\Big(\frac{d^2f}{dt^2}\Big)(t)
\,.
\end{equation}
We also note that
\begin{equation}\label{eq:der-fracI1}
\partial_t^{\alpha}f(t)
=\frac{d}{dt}I^{2-\alpha}\big(f'-f'(0)\big)(t)
\,,
\end{equation}
when $f'$ is absolutely continuous.

For arbitrary constants $\alpha,\beta> 0$, we denote the Mittag--Leffler functions by
\begin{equation}
E_{\alpha,\beta}(z):= \sum_{k=0}^\infty\frac{z^k}{\Gamma(\alpha k+\beta)},
\quad z\in\C.
\end{equation}
One notices that the power series $E_{\alpha,\beta}(z)$ is an entire function
of $z\in\C$. The Mittag--Leffler function $E_{\alpha,1}(z)$ is usually denoted by $E_{\alpha}(z)$. We observe that $E_{\alpha}(0)=1$.

The proof of the following result can be found in \cite[p. 35]{Pod}, see also \cite[Lemma 3.1]{SY}.
In the following we denote the Laplace transform of a function $f(t)$ by the symbol
\begin{equation}\label{eq:laplace}
{\cal L}[f(t)](z):=\int_0^{\infty}e^{-zt}f(t)\ dt
\qquad z\in\C.
\end{equation}
\begin{lemma}
\begin{enumerate}
\item 
Let $\alpha\in(1,2)$ and $\beta>0$ be. Then for any $\mu\in\R$ such that $\pi\alpha/2<\mu<\pi$ there exists a constant $C=C(\alpha,\beta,\mu)>0$ such that 
\begin{equation}\label{eq:stimeE}
\big|E_{\alpha,\beta}(z)\big|\le \frac{C}{1+|z|},
\qquad z\in\C,\ \mu\le|\arg(z)|\le\pi.
\end{equation}
\item
For $\alpha\,,\beta\,,\lambda>0$ one has
\begin{equation}\label{eq:LTML}
{\cal L}\big[t^{\beta-1}E_{\alpha,\beta}(-\lambda t^\alpha)\big)\big](z)=\frac{z^{\alpha-\beta}}{z^\alpha+\lambda}
\qquad \Re z>\lambda^{\frac1\alpha}\,.
\end{equation}
\item 
If $\alpha\,,\lambda>0$, then we have
\begin{equation}\label{eq:Ea1}
\frac{d}{dt}E_{\alpha}(-\lambda t^{\alpha})=-\lambda t^{\alpha-1}E_{\alpha,\alpha}(-\lambda t^{\alpha}), 
\qquad t>0,
\end{equation}
\begin{equation}\label{eq:Eaa1}
\frac{d}{dt}\Big(t^{k}E_{\alpha,k+1}(-\lambda t^{\alpha})\Big)=t^{k-1}E_{\alpha,k}(-\lambda t^{\alpha}), 
\qquad k\in\N\,,t\ge0,
\end{equation}
\begin{equation}\label{eq:Eaaa}
\frac{d}{dt}\Big(t^{\alpha-1}E_{\alpha,\alpha}(-\lambda t^{\alpha})\Big)=t^{\alpha-2}E_{\alpha,\alpha-1}(-\lambda t^{\alpha}), 
\qquad t\ge0.
\end{equation}
\end{enumerate}
\end{lemma}
We recall an elementary result that is useful in the estimates.
\begin{lemma} For any $0<\beta<1$ the function $x\to\frac{x^\beta}{1+x}$ gains its maximum on $[0,+\infty[$ at point $\frac\beta{1-\beta}$ and the maximum value is given by
\begin{equation}\label{eq:maxbeta}
\max_{x\ge0}\frac{x^\beta}{1+x}=\beta^\beta(1-\beta)^{1-\beta},
\qquad\beta\in (0,1)
\,.
\end{equation}

\end{lemma}
We now call to mind the definition of fractional vector-valued Sobolev spaces.
For $\beta\in (0,1)$, $T>0$ and a Hilbert space $H$, endowed with the norm $\|\cdot\|_H$, $H^\beta(0,T;H)$ is the space of all $u\in L^2(0,T;H)$ such that 
\begin{equation}\label{eq:gagliardo}
[u]_{H^\beta(0,T;H)}:=\left(\int_0^T\int_0^T\frac{\|u(t)-u(\tau)\|_H^2}{|t-\tau|^{1+2\beta}}\ dtd\tau\right)^{1/2}<+\infty\,,
\end{equation}
that is $[u]_{H^\beta(0,T;H)}$ is the so-called Gagliardo semi-norm of $u$.
$H^\beta(0,T;H)$ is endowed with the norm
\begin{equation}\label{eq:defHs}
\|\cdot\|_{H^\beta(0,T;H)}:=\|\cdot\|_{L^2(0,T;H)}+[\ \cdot\ ]_{H^\beta(0,T;H)}.
\end{equation}
The following extension of a known result (see
\cite[Theorem 2.1]{GLY}) to the case of vector valued functions is crucial in the proof of Theorem \ref{th:hidstr}.
We use the symbol $\sim$ between norms to indicate two equivalent norms.
\begin{theorem}\label{th:R-Lop}
Let $H$ be a separable Hilbert space.
\begin{itemize}
\item[(i)] 
The Riemann--Liouville operator $I^{\beta}:L^2(0,T;H)\to L^2(0,T;H)$, $0<\beta\le1$, is injective and the range ${\cal R}(I^{\beta})$ of $I^{\beta}$ is given by
\begin{equation}\label{eq:range}
{\cal R}(I^{\beta})=
\begin{cases}H^\beta(0,T;H), \hskip5.5cm 0<\beta<\frac12,
\\
\Big\{v\in H^{\frac12}(0,T;H): \int_0^Tt^{-1}|v(t)|^2 dt<\infty\Big\},
\qquad \beta=\frac12,
\\
_0H^\beta(0,T;H), \hskip5.2cm \frac12<\beta\le1,
\end{cases}
\end{equation}
where $_0H^\beta(0,T)=\{u\in H^\beta(0,T): u(0)=0\}$.
\item[(ii)] 
For the Riemann--Liouville operator $I^{\beta}$ and its inverse operator $I^{-\beta}$ the norm equivalences 
\begin{equation}\label{eq:R-Lop}
\begin{split}
\|I^{\beta}(u)\|_{H^\beta(0,T;H)}
&\sim\|u\|_{L^2(0,T;H)},
\qquad u\in L^2(0,T;H),
\\
\|I^{-\beta}(v)\|_{L^2(0,T;H)}
&\sim\|v\|_{H^\beta(0,T;H)},
\qquad v\in {\cal R}(I^{\beta}),
\end{split}
\end{equation}
hold true.
\end{itemize}
\end{theorem}

\section{Existence and regularity of solutions}\label{s:reg}

First, we recall the definition of  $H^2$-solutions and  strong solutions,  see \cite{LS4}. Let $\alpha\in(1,2)$ and $T>0$.
\begin{definition}\label{de:wss}
\begin{enumerate}
\item[\bf 1.] 
A function $u$ is called a $H^2$-solution of the fractional boundary value problem
\begin{equation}\label{eq:weakp}
\begin{cases}
\partial_t^{\alpha}u+\Delta^2 u=0
\qquad \mbox{in}\ (0,T)\times\Omega,
\\
u=\Delta u=0    \hskip1.4cm \mbox{on}\ (0,T)\times\partial\Omega,
\end{cases}
\end{equation}
if  $u\in C([0,T];H^2(\Omega)\cap H^1_0(\Omega))$, $u_t\in L^2(0,T;L^2(\Omega))\cap C([0,T];D(A^{-\theta}))$,  for some $\theta\in(0,1)$, 
and for any $v\in H^2(\Omega)\cap H^1_0(\Omega)$ one has $\int_\Omega I^{2-\alpha}\big(u_t(\cdot,x)-u_t(0,x)\big)v(x) \ dx\in C^1([0,T])$ and 
\begin{equation}\label{eq:w-int}
\frac{d}{dt}\int_\Omega I^{2-\alpha}\big(u_t(t,x)-u_t(0,x)\big)v(x) \ dx+\int_\Omega\Delta u(t,x)\Delta v(x) \ dx=0
\qquad t\in [0,T]
\,.
\end{equation}
\item [\bf 2.]
A function $u$ is called a strong solution if $u\in C([0,T];D(A))\cap C^1([0,T];L^2(\Omega))$, $\partial_t^{\alpha}u\in C([0,T];L^2(\Omega))$ and satisfies
problem \eqref{eq:weakp}.
\end{enumerate}
\end{definition}
\begin{remark}
A strong solution is also a $H^2$-solution.
\end{remark}
We recall the following existence result, see \cite[Theorem 4.3]{LS4}.
\begin{theorem}\label{th:reg-l2}
\begin{enumerate}
\item[(i)] 
If $u_0\in H^2(\Omega)\cap H^1_0(\Omega)$ and $u_1\in L^{2}(\Omega)$, then the function
\begin{equation}\label{eq:def-u0}
u(t,x)=\sum_{n=1}^\infty\big[  \langle u_0,e_n\rangle E_{\alpha}(-\lambda_nt^\alpha) 
+\langle u_1,e_n\rangle t E_{\alpha,2}(-\lambda_nt^\alpha)\big]e_n(x)
\end{equation}
is the unique $H^2$-solution of \eqref{eq:weakp} satisfying
 the initial conditions
\begin{equation}\label{eq:ini-cond}
u(0,\cdot)=u_{0},\quad
u_t(0,\cdot)=u_{1}.
\end{equation}
In addition
\begin{equation}\label{eq:def-u_t}
u_t(t,x)=\sum_{n=1}^\infty\big[-\lambda_n \langle u_0,e_n\rangle t^{\alpha-1}E_{\alpha,\alpha}(-\lambda_nt^\alpha)
+\langle u_1,e_n\rangle E_{\alpha}(-\lambda_nt^\alpha)\big]e_n(x)
\,,
\end{equation}
and $u_t\in C([0,T];D(A^{-\theta}))$ for  $\theta\in\big(\frac{2-\alpha}{2\alpha},\frac12\big)$.
\item[(ii)]
For $u_0\in H^4(\Omega)\cap \{u\in H^3(\Omega): u=\Delta u=0\ \mbox{on}\ \partial\Omega\}$ and $u_1\in H^2(\Omega)\cap H^1_0(\Omega)$
the $H^2$-solution given by \eqref{eq:def-u0} is a strong one and 
\begin{equation}\label{eq:def-u-alpha}
\partial_t^{\alpha}u(t,x)=-\sum_{n=1}^\infty\big[\lambda_n \langle u_0,e_n\rangle E_{\alpha,1}(-\lambda_nt^\alpha)
+\lambda_n\langle u_1,e_n\rangle t E_{\alpha,2}(-\lambda_nt^\alpha)\big]e_n(x)
\,.
\end{equation}
\end{enumerate}
\end{theorem}
\begin{proposition}
For $u_0\in H^2(\Omega)\cap H^1_0(\Omega)$ and $u_1\in L^{2}(\Omega)$ if $u$ is the $H^2$-solution given by \eqref{eq:def-u0}, then $A^{-\frac12} u$ is a strong solution.
\end{proposition}
\begin{Proof}
To prove the statement we set $w=A^{-\frac12} u$. So, we have $w\in D(A^{\frac12})=H^2(\Omega)\cap H^1_0(\Omega)$, and in particular $w=0$ on $\partial\Omega$. Moreover 
\begin{equation*}
-\triangle w=A^{\frac12}A^{-\frac12} u=u \quad\text{on}\ \Omega.
\end{equation*}
Since $u\in H^2(\Omega)\cap H^1_0(\Omega)$, we get $w\in H^4(\Omega)$ and $\triangle w=0$ on $\partial\Omega$.

In addition, thanks to Theorem \ref{th:reg-l2} we have $A^{\frac12}w=u\in C^1([0,T];D(A^{-\theta}))$ for  $\theta\in\big(\frac{2-\alpha}{2\alpha},\frac12\big)$. Therefore
\begin{equation*}
\|w_t(t,\cdot)\|_{L^2(\Omega)}^2=
\sum_{n=1}^\infty\lambda_n^{2\theta-1}\lambda_n^{1-2\theta} |\langle w_t,e_n\rangle|^2
\le C\sum_{n=1}^\infty\lambda_n^{1-2\theta} |\langle w_t,e_n\rangle|^2
= C \|A^{\frac12}w_t(t,\cdot)\|_{D(A^{-\theta})}^2,
\end{equation*}
and hence  $A^{-\frac12} u=w\in C^1([0,T];L^2(\Omega))$.

To complete the proof we have to show that 
\begin{equation}\label{eq:abs-c}
\frac{d}{dt} I^{2-\alpha}\big(w_t-A^{-\frac12}u_1\big)=\partial_t^{\alpha}w
\,.
\end{equation}
To this end, thanks to \eqref{eq:der-fracI1}  it is sufficient to prove that
$w_t$ is absolutely continuous. Indeed, by \eqref{eq:def-u_t} we get
\begin{equation*}
w_t(t,\cdot)=A^{-\frac12} u_t(t,\cdot)
=\sum_{n=1}^\infty\lambda_n^{-\frac12}\big[-\lambda_n \langle u_0,e_n\rangle t^{\alpha-1}E_{\alpha,\alpha}(-\lambda_nt^\alpha)
+\langle u_1,e_n\rangle E_{\alpha}(-\lambda_nt^\alpha)\big]e_n
\,.
\end{equation*}
The series of the derivatives with respect to the variable $t$, thanks to \eqref{eq:Eaaa} and \eqref{eq:Ea1}, is given by
\begin{equation}\label{eq:un''}
-\sum_{n=1}^\infty\lambda_n^{\frac12}\big[ \langle u_0,e_n\rangle t^{\alpha-2}E_{\alpha,\alpha-1}(-\lambda_nt^\alpha)
+\langle u_1,e_n\rangle t^{\alpha-1} E_{\alpha,\alpha}(-\lambda_nt^\alpha)\big]e_n.
\end{equation}
In virtue of \eqref{eq:stimeE} we note that
\begin{multline*}
\int_0^T\big\|\sum_{n=1}^\infty \lambda_n^{\frac12}\big[ \langle u_0,e_n\rangle t^{\alpha-2}E_{\alpha,\alpha-1}(-\lambda_nt^\alpha)
+\langle u_1,e_n\rangle t^{\alpha-1} E_{\alpha,\alpha}(-\lambda_nt^\alpha)\big]e_n\big\|_{L^2(\Omega)}\ dt
\\
\le
C(T^{\alpha-1}\|A^{\frac12} u_0\|_{L^2(\Omega)}+T^{\frac\alpha2}\| u_1\|_{L^2(\Omega)}),
\end{multline*}
that is the series in \eqref{eq:un''} belongs to $L^1(0,T;L^2(\Omega))$,  and hence $w_t$ is absolutely continuous.

In conclusion,  if we substitute in formula \eqref{eq:w-int}  $u$ with $A^{\frac12}w$, thanks to \eqref{eq:abs-c} we obtain
\begin{equation*}
\int_\Omega A^{\frac12}\big(\partial_t^{\alpha}w(t,x)+\Delta^2 w(t,x)\big)v(x) \ dx=0
\qquad t\in [0,T], \ \  v\in H^2(\Omega)\cap H^1_0(\Omega)
\,,
\end{equation*}
that is $\partial_t^{\alpha}w+\Delta^2 w=0$. 
\end{Proof}
\begin{definition}\label{de:wss1}
\begin{enumerate}
\item[\bf 1.] 
A function $u$ is called a $H^3$-solution of the fractional boundary value problem
\begin{equation}\label{eq:weakp1}
\begin{cases}
\partial_t^{\alpha}u+\Delta^2 u=0
\qquad \mbox{in}\ (0,T)\times\Omega,
\\
u=\Delta u=0    \hskip1.4cm \mbox{on}\ (0,T)\times\partial\Omega,
\end{cases}
\end{equation}
if  
$u\in C([0,T];D(A^{\frac34}))$, $u_t\in L^2(0,T;H^1(\Omega))\cap C([0,T];H^{-1}(\Omega))$,  $\partial_t^{\alpha} u\in L^2([0,T];L^2(\Omega))$
and for any $v\in H^1_0(\Omega)$ one has 
\begin{equation}
\int_\Omega \partial_t^{\alpha} u(t,x)v(x) \ dx-\int_\Omega\nabla\Delta u(t,x)\cdot\nabla v(x) \ dx=0 
\qquad t\in (0,T)
\,.
\end{equation}
\item [\bf 2.]
A function $u\in C([0,T];H_0^1(\Omega))$ is called a $H^1$-solution of \eqref{eq:weakp1} if $A^{-\frac12} u$ is a $H^3$-solution.
\end{enumerate}
\end{definition}
\begin{theorem}\label{th:reg-l21}
\begin{enumerate}
\item[(i)] 
Let $u_0\in D(A^{\frac34}) $ and $u_1\in H^1_0(\Omega)$. Then the $H^2$-solution 
\begin{equation}
u(t,x)=\sum_{n=1}^\infty\big[  \langle u_0,e_n\rangle E_{\alpha}(-\lambda_nt^\alpha) 
+\langle u_1,e_n\rangle t E_{\alpha,2}(-\lambda_nt^\alpha)\big]e_n(x)
\end{equation}
of \eqref{eq:weakp1} with initial conditions
\begin{equation}\label{Aeq:cauchy10}
u(0,\cdot)=u_{0},\quad
u_t(0,\cdot)=u_{1},
\end{equation}
is a  $H^3$-solution  
and
\begin{equation}\label{eq:partial-l2}
\|\partial_t^{\alpha}u\|_{L^2(0,T;L^2(\Omega))}\le C\big(\|\nabla\Delta u_0\|_{L^2(\Omega)}+\|\nabla u_1\|_{L^2(\Omega)}\big)\,,
\qquad (C>0).
\end{equation}
In addition, for $\theta\in\big(0,\frac1{2\alpha}\big)$ we have
$\nabla\Delta u\in L^2(0,T;D(A^{\theta}))$ and
\begin{equation}\label{eq:nabla-l2}
\|\nabla\Delta u\|_{L^2(0,T;D(A^{\theta}))}\le C\big(\|\nabla\Delta u_0\|_{L^2(\Omega)}+\|\nabla u_1\|_{L^2(\Omega)}\big)\,,
\qquad (C>0).
\end{equation}
\item[(ii)] 
Let $u_0\in H^1_0(\Omega)$ and $u_1\in H^{-1}(\Omega)$. Then the function
\begin{equation}
u(t,x)=\sum_{n=1}^\infty\big[  \langle u_0,e_n\rangle E_{\alpha}(-\lambda_nt^\alpha) 
+\langle u_1,e_n\rangle_{-\frac14,\frac14} t E_{\alpha,2}(-\lambda_nt^\alpha)\big]e_n(x)
\end{equation}
is a  $H^1$-solution  of \eqref{eq:weakp1}.
\end{enumerate}
\end{theorem}

\begin{Proof}
(i) First, we note that for any $t\in [0,T]$ we have $ u(t)\in D(A^{\frac34})$. Indeed, in view of \eqref{eq:H^{-1}0} we have
\begin{equation*}
\|\nabla\Delta  u(t)\|_{L^2(\Omega)}^2
\le
2\sum_{n=1}^\infty \lambda_n^{\frac32} \big| \langle u_0,e_n\rangle 
E_{\alpha}(-\lambda_nt^\alpha) \big|^2
+2\sum_{n=1}^\infty \lambda_n^{\frac32} \big| \langle u_1,e_n\rangle t 
E_{\alpha,2}(-\lambda_nt^\alpha)\big|^2.
\end{equation*}
Thanks to \eqref{eq:stimeE} we get
\begin{equation*}
\begin{split}
 \lambda_n^{\frac32} \big| \langle u_0,e_n\rangle E_{\alpha}(-\lambda_nt^\alpha) \big|^2
 &\le C\lambda_n^{\frac32} \big| \langle u_0,e_n\rangle  \big|^2,
\\
 \lambda_n^{\frac32} \big|\langle u_1,e_n\rangle t E_{\alpha,2}(-\lambda_nt^\alpha)\big|^2
& \le C t^{2-\alpha} \lambda_n^{\frac12}\big|\langle u_1,e_n\rangle \big|^2\frac{\lambda_nt^\alpha}{(1+\lambda_nt^\alpha)^2}
\le C t^{2-\alpha} \lambda_n^{\frac12}\big|\langle u_1,e_n\rangle \big|^2,
 \end{split}
\end{equation*}
whence, being $\alpha<2$, we obtain
\begin{equation}\label{eq:deltau}
\|\nabla\Delta u(t)\|_{L^2(\Omega)}^2
\le C\|\nabla\Delta u_0\|_{L^2(\Omega)}^2+C T^{2-\alpha} \|\nabla u_1\|_{L^2(\Omega)}^2.
\end{equation}
Following the same argumentations used to prove \eqref{eq:deltau}, we get for any $n\in\N$ 
\begin{multline*}
\Big\|\nabla\Delta \sum_{k=n}^\infty \big[  \langle u_0,e_k\rangle E_{\alpha}(-\lambda_kt^\alpha) 
+\langle u_1,e_k\rangle t E_{\alpha,2}(-\lambda_kt^\alpha)\big]e_k\Big\|_{L^2(\Omega)}^2
\\
\le C\sum_{k=n}^\infty \lambda_k^{\frac32} \big| \langle u_0,e_k\rangle  \big|^2+C T^{2-\alpha} 
\sum_{k=n}^\infty \lambda_k^{\frac12}\big|\langle u_1,e_k\rangle \big|^2,
\end{multline*}
and hence
\begin{equation*}
\lim_{n\to\infty}\sup_{t\in [0,T]}
\Big\|\nabla\Delta \sum_{k=n}^\infty \big[  \langle u_0,e_k\rangle E_{\alpha}(-\lambda_kt^\alpha) 
+\langle u_1,e_k\rangle t E_{\alpha,2}(-\lambda_kt^\alpha)\big]e_k\Big\|_{L^2(\Omega)}=0\,.
\end{equation*}
As a consequence, the series 
$\sum_{n=1}^\infty\big[  \langle u_0,e_n\rangle E_{\alpha}(-\lambda_nt^\alpha) 
+\langle u_1,e_n\rangle t E_{\alpha,2}(-\lambda_nt^\alpha)\big]e_n$ 
is convergent in $D(A^{\frac34})$ uniformly in $t\in [0,T]$, so $u\in C([0,T];D(A^{\frac34}))$.
\begin{multline*}
\Big\|\frac{d}{dt}I^{2-\alpha}\big(u_t(t,\cdot)-u_1\big)\Big\|_{H^{-1}(\Omega)}^2=
\sum_{n=1}^\infty\lambda_n ^{-\frac12}\Big|\frac{d}{dt}\Big(\lambda_n \langle u_0,e_n\rangle t E_{\alpha,2}(-\lambda_nt^\alpha)
+\lambda_n \langle u_1,e_n\rangle t^2 E_{\alpha,3}(-\lambda_nt^\alpha)\Big)\Big|^2
\\
\le C\big(\|\nabla\Delta u_0\|_{L^2(\Omega)}+t^{2-\alpha}\|\nabla u_1\|_{L^2(\Omega)}\big)
\,.
\end{multline*}

To prove \eqref{eq:partial-l2} we note that,
thanks to \eqref{eq:def-u-alpha} and \eqref{eq:stimeE}, we have
\begin{multline}\label{eq:par-u-alpha}
\|\partial_t^{\alpha}u(t,\cdot)\|_{L^2(\Omega)}^2
=\sum_{n=1}^\infty\big|\lambda_n \langle u_0,e_n\rangle E_{\alpha}(-\lambda_nt^\alpha)
+\lambda_n\langle u_1,e_n\rangle t E_{\alpha,2}(-\lambda_nt^\alpha)\big|^2
\\
\le
C\sum_{n=1}^\infty\lambda_n^{\frac32} |\langle u_0,e_n\rangle|^2 \frac{\lambda_n^{\frac12}}{(1+\lambda_nt^\alpha)^2}
+C\sum_{n=1}^\infty \lambda_n^{\frac12}|\langle u_1,e_n\rangle|^2 \frac{\lambda_n^{\frac32}t^2}{(1+\lambda_nt^\alpha)^2}\,.
\end{multline}
Observing that
\begin{equation*}
\frac{\lambda_n^{\frac12}}{(1+\lambda_nt^\alpha)^2}=\Big(\frac{(\lambda_nt^\alpha)^{\frac{1}4}}{1+\lambda_nt^\alpha}\Big)^2 t^{-\frac\alpha2},
\qquad\frac{\lambda_n^{\frac32}t^2}{(1+\lambda_nt^\alpha)^2}
=\Big(\frac{(\lambda_nt^\alpha)^{\frac34}}{1+\lambda_nt^\alpha}\Big)^2 t^{2-\frac32\alpha},
\end{equation*}
from  \eqref{eq:maxbeta} and \eqref{eq:par-u-alpha} we deduce
\begin{equation*}
\|\partial_t^{\alpha}u(\cdot,t)\|_{L^2(\Omega)}^2
\le C t^{-\frac\alpha2}\ \|\nabla \Delta u_0\|_{L^2(\Omega)}^2+Ct^{2-\frac32\alpha}\|\nabla u_1\|_{L^2(\Omega)}^2\,.
\end{equation*}
Therefore, since $\alpha<2$ we have 
\begin{equation*}
\|\partial_t^{\alpha}u\|_{L^2(0,T;L^2(\Omega))}^2
\le C T^{1-\frac\alpha2}\ \|\nabla \Delta u_0\|_{L^2(\Omega)}^2+CT^{3-\frac32\alpha}\|\nabla u_1\|_{L^2(\Omega)}^2\,,
\end{equation*}
so
$\partial_t^{\alpha} u\in L^2(0,T;L^2(\Omega))$
and  \eqref{eq:partial-l2} holds.

Now, we fix $\theta\in\big(0,\frac1{2\alpha}\big)$. Thanks to \eqref{eq:def-u0} and \eqref{eq:stimeE} we get
\begin{multline*}
\|\nabla\Delta u(t,\cdot)\|_{D(A^{\theta})}^2=\sum_{n=1}^\infty\lambda_n^{\frac32+2\theta}\big| \langle u_0,e_n\rangle E_{\alpha}(-\lambda_nt^\alpha)
+\langle u_1,e_n\rangle t E_{\alpha,2}(-\lambda_nt^\alpha)\big|^2
\\
\le
C\sum_{n=1}^\infty\lambda_n^{\frac32} |\langle u_0,e_n\rangle|^2 \frac{\lambda_n^{2\theta}}{(1+\lambda_nt^\alpha)^2}
+C\sum_{n=1}^\infty\lambda_n^{\frac12}|\langle u_1,e_n\rangle|^2 \frac{\lambda_n^{1+2\theta}t^2}{(1+\lambda_nt^\alpha)^2}
\,.
\end{multline*}
Since
\begin{equation*}
 \frac{\lambda_n^{2\theta}}{(1+\lambda_nt^\alpha)^2}=\Big(\frac{(\lambda_nt^\alpha)^{\theta}}{1+\lambda_nt^\alpha}\Big)^2 t^{-2\alpha\theta},
\qquad
\frac{\lambda_n^{1+2\theta}t^2}{(1+\lambda_nt^\alpha)^2}=\Big(\frac{(\lambda_nt^\alpha)^{\frac{1+2\theta}2}}{1+\lambda_nt^\alpha}\Big)^2 t^{2-\alpha(1+2\theta)},
\end{equation*}
and $0<\theta<\frac12$, we  can apply \eqref{eq:maxbeta} to have
\begin{equation*}
\|\nabla\Delta u(t,\cdot)\|_{D(A^{\theta})}^2
\le
Ct^{-2\alpha\theta}\|\nabla\Delta u_0\|_{L^2(\Omega)}^2
+Ct^{2-\alpha(1+2\theta)}\|\nabla u_1\|_{L^2(\Omega)}^2.
\end{equation*}
Taking into account that $\theta\in\big(0,\frac1{2\alpha}\big)$ we have $\nabla\Delta u\in L^2(0,T;D(A^{\theta}))$
and  \eqref{eq:nabla-l2} follows.

(ii)

\end{Proof}

\section{Hidden regularity results}\label{s:hidreg}
We develop a procedure similar to that  followed in \cite{K} for Petrovsky systems.
First, we prove some identities that are useful in the proof of the main theorem. 
\begin{lemma}\label{le:tech0}
If $w\in H^4(\Omega)$ with $\Delta w=0$ on $\partial\Omega$
and $h:\overline{\Omega}\to\R^N$ is a vector field of class $C^1$, then 
\begin{multline}\label{eq:triangleuF00}
2\int_\Omega
\Delta^2 w\
h\cdot\nabla\Delta w\ dx 
=
\int_{\partial\Omega}h\cdot\nu |\partial_\nu \Delta w|^2\ d\sigma 
-2
\sum_{i,j=1}^N\int_\Omega
\partial_i  h_j\partial_i \Delta w\partial_j\Delta w\ dx 
\\
+\int_{\Omega}
\sum_{j=1}^N \partial_jh_j\ |\nabla\Delta w|^2\ dx 
\,.
\end{multline}
\end{lemma}
\begin{Proof}
First, we observe that by $w\in H^4(\Omega)$ it follows that the normal derivative $\partial_\nu\Delta w$ is well defined on $\partial\Omega$.

Integrating by parts we get
\begin{equation*}
\int_\Omega
\Delta^2 w\
h\cdot\nabla\Delta w\ dx
=
\int_{\partial\Omega}
\partial_\nu\Delta w\
h\cdot\nabla\Delta w\ d\sigma
-\int_\Omega\nabla\Delta w
\cdot\nabla \big( h\cdot\nabla\Delta w\big)\ dx
\,.
\end{equation*}
Since $\Delta w=0$ on $\partial\Omega$ we have
\begin{equation*}
\nabla\Delta w=(\partial_\nu\Delta w)\nu
\quad
\mbox{on}\ \partial\Omega,
\end{equation*}
(see e.g. \cite[Lemma 2.1]{MM} for a detailed proof)
and hence
\begin{equation}\label{eq:triangleu00}
\int_\Omega
\Delta^2 w\
h\cdot\nabla\Delta w\ dx
=
\int_{\partial\Omega}h\cdot\nu |\partial_\nu \Delta w|^2\ d\sigma
-\int_\Omega\nabla\Delta w
\cdot\nabla \big( h\cdot\nabla\Delta w\big)\ dx
\,.
\end{equation}
We note that
\begin{multline*}
\int_\Omega\nabla\Delta w\cdot\nabla \big( h\cdot\nabla\Delta w \big)\ dx
\\
=\sum_{i,j=1}^N\int_\Omega
\partial_i \Delta w\ \partial_i ( h_j\partial_j\Delta w)\ dx
=\sum_{i,j=1}^N\int_\Omega\partial_i\Delta w\ \partial_i  h_j\partial_j\Delta w\ dx
+\sum_{i,j=1}^N\int_\Omega h_j\ \partial_i\Delta w\partial_j ( \partial_i\Delta w)\ dx.
\end{multline*}
We estimate the last term on the right-hand side again by an integration  by parts, so we obtain
\begin{equation*}
\begin{split}
\sum_{i,j=1}^N\int_\Omega
h_j\ \partial_i\Delta w \partial_j ( \partial_i\Delta w)\ dx
=&\frac12\sum_{j=1}^N\int_\Omega
h_j\ \partial_j \Big( \sum_{i=1}^N(\partial_i\Delta w)^2\Big)\ dx
\\
=&\frac12\int_{\partial\Omega}
h\cdot\nu |\nabla\Delta w|^2\ d\sigma
-\frac12\int_{\Omega}
\sum_{j=1}^N \partial_jh_j\ |\nabla\Delta w|^2\ dx
\,.
\end{split}
\end{equation*}
In conclusion, putting the above identities into \eqref{eq:triangleu00} we deduce \eqref{eq:triangleuF00}.
\end{Proof}

We prove the next result for strong solutions, see Definition \ref{de:wss}-2, whose existence is given by Theorem \ref{th:reg-l2}--(ii). We recall that $I^{\beta}$ is the Riemann--Liouville operator of order $\beta>0$, see \eqref{eq:RLfi}.
\begin{lemma}\label{le:tech}
Suppose $u$ is a strong solution  of 
\begin{equation}\label{eq:stato}
\begin{cases}
\partial_t^{\alpha}u+\Delta^2 u=0
\qquad \mbox{in}\ [0,T]\times\Omega,
\\
u=\Delta u=0    \hskip1.4cm \mbox{on}\ [0,T]\times\partial\Omega.
\end{cases}
\end{equation}
For a vector field $h:\overline{\Omega}\to\R^N$ of class $C^1$ and $\beta\in(0,1)$ we have
\begin{multline}\label{eq:identity}
\int_{\partial\Omega}
h\cdot\nu \big|  I^{\beta}(\partial_\nu\Delta u)(t)\big|^2
\ d\sigma
=-2\int_{\Omega} I^{\beta}(\partial_t^{\alpha}u)(t) h\cdot I^{\beta}(\nabla\Delta u)(t)\ dx
\\
+2\sum_{i,j=1}^N\int_\Omega
\partial_i  h_j I^{\beta}(\partial_i\Delta u)(t)I^{\beta}(\partial_j \Delta u)(t)\ dx
-\int_{\Omega}
\sum_{j=1}^N \partial_jh_j\ | I^{\beta}(\nabla\Delta u)(t)|^2\ dx
\,, \qquad t\in [0,T],
\end{multline}
\begin{multline}\label{eq:identity2}
\int_{\partial\Omega}h\cdot\nu \big| I^{\beta}(\partial_\nu\Delta u )(t)-I^{\beta}(\partial_\nu\Delta u )(\tau)\big|^2 d\sigma 
\\
=
-2\int_{\Omega}\big( I^{\beta}(\partial_t^{\alpha}u)(t)-I^{\beta}(\partial_{t}^{\alpha}u)(\tau)\big) h\cdot \big(I^{\beta}(\nabla\Delta u)(t)-I^{\beta}(\nabla\Delta u)(\tau)\big)\ dx
\\
+2
\sum_{i,j=1}^N\int_\Omega
\partial_i  h_j  \big(I^{\beta}(\partial_i\Delta u)(t)-I^{\beta}(\partial_i\Delta u)(\tau)\big)  \big(I^{\beta}(\partial_j \Delta u)(t)-I^{\beta}(\partial_j\Delta u)(\tau)\big)\ dx 
\\
-
\sum_{j=1}^N\int_{\Omega} \partial_jh_j\ |I^{\beta}(\nabla\Delta u )(t)-I^{\beta}(\nabla\Delta u )(\tau)|^2\ dx 
\,, \qquad t,\tau\in [0,T]
\,.
\end{multline}
\end{lemma}
\begin{Proof}
First, we apply the operator $I^{\beta}$, $\beta\in(0,1)$, to  equation \eqref{eq:stato}:
\begin{equation}\label{eq:eqI}
I^{\beta}(\partial_t^{\alpha}u)(t)=-I^{\beta}(\Delta^2 u )(t)
\qquad t\in [0,T].
\end{equation}
By means of the scalar product in $L^2(\Omega)$ 
we multiply  \eqref{eq:eqI} by 
\begin{equation*}
2h\cdot \nabla\Delta I^{\beta}( u)(t),
\end{equation*}
that is
\begin{equation}\label{eq:identity1}
2\int_\Omega I^{\beta}(\partial_t^{\alpha}u)(t)h\cdot\nabla\Delta I^{\beta}(u)(t)\ dx
=
-2\int_\Omega \Delta^2 I^{\beta}( u )(t)
h\cdot\nabla\Delta I^{\beta}(u)(t)\ dx  
\,.
\end{equation}
To evaluate the term
\begin{equation*}
2\int_\Omega\Delta^2 I^{\beta}( u )(t)
h\cdot\nabla\Delta I^{\beta}(u)(t)\ dx  
\,,
\end{equation*} 
we apply Lemma \ref{le:tech0} to the function
$
w(t,x)=I^{\beta}( u )(t)
$,
so
from \eqref{eq:triangleuF00} we deduce
\begin{multline*}
2\int_\Omega
\Delta^2 I^{\beta}( u )(t)
h\cdot\nabla\Delta I^{\beta}(u)(t)\ dx  
=\int_{\partial\Omega}h\cdot\nu \big| I^{\beta}(\nabla\Delta u )(t)\big|^2 d\sigma 
\\
-2
\sum_{i,j=1}^N\int_\Omega
\partial_i  h_j  I^{\beta}(\partial_i\Delta u)(t)  I^{\beta}(\partial_j \Delta u)(t)\ dx 
+
\int_{\Omega} \sum_{j=1}^N\partial_jh_j\ |I^{\beta}(\nabla\Delta u )(t)|^2\ dx 
\,.
\end{multline*}
In conclusion, plugging the above formula into \eqref{eq:identity1}, we obtain \eqref{eq:identity}.

The proof of \eqref{eq:identity2} is alike: we start from 
\begin{equation*}
I^{\beta}(\partial_t^{\alpha}u)(t)-I^{\beta}(\partial_{t}^{\alpha}u)(\tau)=-\big(I^{\beta}(\Delta^2 u )(t)-I^{\beta}(\Delta^2 u )(\tau)\big)
\qquad \qquad t,\tau\in [0,T],
\end{equation*}
and we multiply both terms  by
$
2h\cdot\nabla\Delta \big(I^{\beta}(u)(t)-I^{\beta}(u)(\tau)\big).
$
Then, applying  Lemma \ref{le:tech0} to the function
$
w(t,\tau,x)=I^{\beta}( u )(t)-I^{\beta}( u )(\tau)
$ we get the  identity \eqref{eq:identity2}. 
\end{Proof}

\begin{remark}
We observe that the proof of the identities \eqref{eq:identity} and \eqref{eq:identity2} cannot be done for a general function $w$ and then applied to $w=I^\beta (u)$, since
\begin{equation*}
\partial_t^{\alpha}I^{\beta}(u)\not=I^{\beta}(\partial_t^{\alpha}u),
\end{equation*}
as one can conclude from \eqref{eq:der-frac}.
\end{remark}
\begin{theorem}\label{th:hidstr}
Let  $u_0\in H^4(\Omega)\cap \{u\in H^3(\Omega): u=\Delta u=0\ \mbox{on}\ \partial\Omega\}$ and $u_1\in H^2(\Omega)\cap H^1_0(\Omega)$. Then, for  $T>0$ the strong solution $u$ of problem
\begin{equation}\label{eq:cauchy100}
\begin{cases}
\partial_t^{\alpha}u+\Delta^2 u=0
\qquad \mbox{in}\ [0,T]\times\Omega,
\\
u=\Delta u=0    \hskip1.4cm \mbox{on}\ [0,T]\times\partial\Omega
\\
u(0,\cdot)=u_{0},\quad
u_t(0,\cdot)=u_{1}.
\end{cases}
\end{equation}
 satisfies the inequality
\begin{equation}\label{eq:hidden-alpha}
\int_0^T\int_{\partial\Omega} \big|\partial_\nu\Delta u\big|^2d\sigma dt
\le C\big(\|\nabla\Delta u_0\|^2_{L^2(\Omega)}+\|\nabla u_1\|^2_{L^2(\Omega)}\big)
\,,
\end{equation}
for some constant $C=C(T)>0$.

\end{theorem}

\begin{Proof}
We apply Theorem \ref{th:R-Lop} to prove the statement. Indeed, for $H=L^2(\partial\Omega)$ and $\beta\in(0,1)$
we can apply \eqref{eq:R-Lop}  to have
\begin{equation}
\|\partial_\nu\Delta u\|_{L^2(0,T;L^2(\partial\Omega))}\sim
\|I^{\beta}(\partial_\nu\Delta u)\|_{H^{\beta}(0,T;L^2(\partial\Omega))}\,.
\end{equation}
So, the inequality \eqref{eq:hidden-alpha} is equivalent to  
\begin{equation}\label{eq:hidden-alpha1}
\|I^{\beta}(\partial_\nu\Delta u)\|_{H^{\beta}(0,T;L^2(\partial\Omega))}
\le C\big(\|\nabla\Delta u_0\|_{L^2(\Omega)}+\|\nabla u_1\|_{L^2(\Omega)}\big)
\,.
\end{equation}
To prove \eqref{eq:hidden-alpha1}, taking into account \eqref{eq:defHs},
we  have to estimate
$\big\|I^{\beta}(\partial_\nu\Delta u)\big\|_{L^2(0,T;L^2(\partial\Omega))}$ and 
\break
$\big[ I^{\beta}(\partial_\nu\Delta u) \big]_{H^{\beta}(0,T;L^2(\partial\Omega))}$.
To this end we  employ the two identities in Lemma \ref{le:tech} by means of a suitable choice of the vector field $h$. 
Indeed, we consider a vector field $h\in C^1(\overline{\Omega};\R^N)$ satisfying the condition
\begin{equation}\label{eq:h}
h=\nu
\qquad
\text{on}\quad\partial\Omega
\end{equation}
(see e.g. \cite {K} for the existence of such vector field $h$). If we
integrate  \eqref{eq:identity} over $[0,T]$, then we obtain
\begin{multline*}
\int_0^T\int_{\partial\Omega}
\big|I^{\beta}(\partial_\nu\Delta u)(t)\big|^2
\ d\sigma dt
=-2\int_0^T\int_\Omega I^{\beta}(\partial_t^{\alpha}u)(t) h\cdot I^{\beta}(\nabla\Delta u)(t)dx\ dt
\\
+2\sum_{i,j=1}^N\int_0^T\int_\Omega
\partial_i  h_j I^{\beta}(\partial_i\Delta u)(t)I^{\beta}(\partial_j\Delta u)(t)\ dx dt
-\int_0^T\int_{\Omega}
\sum_{j=1}^N \partial_jh_j\ | I^{\beta}(\nabla\Delta u)(t)|^2\ dxdt
\,.
\end{multline*}
Since $h\in C^1(\overline{\Omega};\R^N)$ from the above inequality we get
\begin{equation}\label{eq:norml2beta}
\big\|I^{\beta}(\partial_\nu\Delta u)\big\|_{L^2(0,T;L^2(\partial\Omega))}
\le C\Big(\big\|I^{\beta}(\partial_t^{\alpha}u)\big\|_{L^2(0,T;L^2(\Omega))}
+\big\|I^{\beta}(\nabla\Delta u)\big\|_{L^2(0,T;L^2(\Omega))}\Big)
\,,
\end{equation}
for some constant $C>0$.

Now we have to evaluate the Gagliardo semi-norm 
$\big[ I^{\beta}(\partial_\nu\Delta u) \big]_{H^{\beta}(0,T;L^2(\partial\Omega))}$, recalled in \eqref{eq:gagliardo}.
To this end we multiply both terms of \eqref{eq:identity2} by $\frac1{|t-\tau|^{1+2\beta}}$ and integrate over $[0,T]\times[0,T]$, then we have
 \begin{multline}\label{eq:identity21}
\big[ I^{\beta}(\partial_\nu\Delta u) \big]_{H^{\beta}(0,T;L^2(\partial\Omega))}^2
\\
=
-2\int_0^T\int_0^T\frac1{|t-\tau|^{1+2\beta}}\int_\Omega\big( I^{\beta}(\partial_t^{\alpha}u)(t)-I^{\beta}(\partial_{t}^{\alpha}u)(\tau)\big)
h\cdot \big(I^{\beta}(\nabla\Delta u)(t)-I^{\beta}(\nabla\Delta u)(\tau)\big)dx\ dt d\tau
\\
+2\int_0^T\int_0^T\frac1{|t-\tau|^{1+2\beta}}
\sum_{i,j=1}^N\int_\Omega
\partial_i  h_j  \big(I^{\beta}(\partial_i\Delta u)(t)-I^{\beta}(\partial_i\Delta u)(\tau)\big)  \big(I^{\beta}(\partial_j\Delta u)(t)-I^{\beta}(\partial_j\Delta u)(\tau)\big)\ dx \ dt d\tau
\\
-
\int_0^T\int_0^T\frac1{|t-\tau|^{1+2\beta}}\int_{\Omega} \sum_{j=1}^N\partial_jh_j\ |I^{\beta}(\nabla\Delta u )(t)-I^{\beta}(\nabla\Delta u )(\tau)|^2\ dx \ dt d\tau
\,.
\end{multline}
The first term on the right-hand side can be estimated as follows
\begin{multline*}
-2\int_0^T\int_0^T\frac1{|t-\tau|^{1+2\beta}}\int_\Omega\big( I^{\beta}(\partial_t^{\alpha}u)(t)-I^{\beta}(\partial_{t}^{\alpha}u)(\tau)\big)
h\cdot \big(I^{\beta}(\nabla\Delta u)(t)-I^{\beta}(\nabla\Delta u)(\tau)\big)dx\ dt d\tau
\\
\le
C\Big(\big[I^{\beta}(\partial_t^{\alpha}u)\big]_{H^{\beta}(0,T;L^2(\Omega))}^2
+\big[I^{\beta}(\nabla\Delta u)\big]_{H^{\beta}(0,T;L^2(\Omega))}^2\Big)
\,,
\end{multline*}
and hence from \eqref{eq:identity21} we deduce
\begin{equation}\label{eq:normhbeta}
\big[ I^{\beta}(\partial_\nu\Delta u) \big]_{H^{\beta}(0,T;L^2(\partial\Omega))}\le
C\Big(\big[I^{\beta}(\partial_t^{\alpha}u)\big]_{H^{\beta}(0,T;L^2(\Omega))}
+\big[I^{\beta}(\nabla\Delta u)\big]_{H^{\beta}(0,T;L^2(\Omega))}\Big)
\,.
\end{equation}
Combining \eqref{eq:norml2beta} and \eqref{eq:normhbeta} we obtain
\begin{equation}\label{eq:normlhbeta}
\big\|I^{\beta}(\partial_\nu\Delta u) \big\|_{H^{\beta}(0,T;L^2(\partial\Omega))}\le
C\Big(\|I^{\beta}(\partial_t^{\alpha}u)\big\|_{H^{\beta}(0,T;L^2(\Omega))}
+\big\|I^{\beta}(\nabla\Delta u)\big\|_{H^{\beta}(0,T;L^2(\Omega))}\Big)
\,.
\end{equation}
Thanks again to Theorem \ref{th:R-Lop} we have
\begin{equation*}
\begin{split}
\big\|I^{\beta}(\partial_t^{\alpha}u)\big\|_{H^{\beta}(0,T;L^2(\Omega))}
&\sim
\|\partial_t^{\alpha}u\|_{L^2(0,T;L^2(\Omega))}
\,,
\\
\big\|I^{\beta}(\nabla\Delta u)\big\|_{H^{\beta}(0,T;L^2(\Omega))}
&\sim
\|\nabla\Delta u\|_{L^2(0,T;L^2(\Omega))}
\,,
\end{split}
\end{equation*}
and hence from \eqref{eq:partial-l2}, \eqref{eq:nabla-l2} and \eqref{eq:normlhbeta}  we deduce
\eqref{eq:hidden-alpha1}. The proof is complete.
\end{Proof}

\begin{theorem}\label{th:hidalpha}
Let  $u_0\in H^1_0(\Omega)$ and $u_1\in H^{-1}(\Omega)$. If $u$ is the $H^1$-solution  of 
\begin{equation}\label{eq:cauchy1}
\begin{cases}
\partial_t^{\alpha}u+\Delta^2 u=0
\qquad \mbox{in}\ (0,T)\times\Omega,
\\
u=\Delta u=0    \hskip1.4cm \mbox{on}\ (0,T)\times\partial\Omega,
\\
u(0,\cdot)=u_{0},\quad
u_t(0,\cdot)=u_{1},
\end{cases}
\end{equation}
then we define the normal derivative $\partial_\nu u$ of $u$ such that  for any $T>0$ we have
\begin{equation}\label{eq:hidden-alpha0}
\int_0^T\int_{\partial\Omega} \big|\partial_\nu u\big|^2d\sigma dt
\le C\Big(\| u_0\|^2_{H^1_0(\Omega)}+\|u_1\|^2_{H^{-1}(\Omega)}\Big)
\,,
\end{equation}
for some constant $C=C(T)$ independent of the initial data.
\end{theorem}
\begin{Proof} 
For $u_0\in H^2(\Omega)\cap H^1_0(\Omega)$ and $u_1\in L^{2}(\Omega)$ we consider the $H^2$- solution $u$ of  
\eqref{eq:cauchy1}. We note that $A^{-\frac12} u$ is the strong solution of  
\eqref{eq:cauchy1} with the initial conditions $u_0$ and $u_1$ replaced by $A^{-\frac12}u_0\in H^4(\Omega)\cap D(A^{\frac34})$ and $A^{-\frac12}u_1\in H^2(\Omega)\cap H^1_0(\Omega)$.
Therefore, by Theorem \ref{th:hidstr} the inequality
\begin{equation*}
\int_0^T\int_{\partial\Omega} \big|\partial_\nu\Delta A^{-\frac12}u\big|^2d\sigma dt
\le C\big(\|\nabla\Delta A^{-\frac12}u_0\|^2_{L^2(\Omega)}+\|\nabla A^{-\frac12}u_1\|^2_{L^2(\Omega)}\big)
\,,
\end{equation*}
holds for any $T>0$. Taking into account that 
\begin{equation*}
\|\nabla A^{-\frac12}u_1\|_{L^2(\Omega)}
=\|u_1\|_{H^{-1}(\Omega)},
\end{equation*}
we get
\begin{equation*}
\int_0^T\int_{\partial\Omega} \big|\partial_\nu u\big|^2d\sigma dt
\le C\Big(\|\nabla u_0\|^2_{L^2(\Omega)}+\|u_1\|^2_{H^{-1}(\Omega)}\Big)
\,.
\end{equation*}
By density there exists a unique continuous linear map
\begin{equation*}
{\cal D_\nu}:  H^1_0(\Omega)\times H^{-1}(\Omega)\to L^2(0,T;L^2(\partial\Omega))
\end{equation*}
such that 
\begin{equation*}
{\cal D_\nu}(u_0,u_1)=\partial_\nu u
\qquad\forall(u_0,u_1)\in H^2(\Omega)\cap H^1_0(\Omega)\times L^2(\Omega)
\end{equation*}
and 
\begin{equation*}
\int_0^T\int_{\partial\Omega} \big|{\cal D_\nu}(u_0,u_1)\big|^2d\sigma dt
\le C\Big(\| u_0\|^2_{H^1_0(\Omega)}+\|u_1\|^2_{H^{-1}(\Omega)}\Big)
\qquad\forall(u_0,u_1)\in H^1_0(\Omega)\times H^{-1}(\Omega)\,.
\end{equation*}
In conclusion, given  $u_0\in H^1_0(\Omega)$ and $u_1\in H^{-1}(\Omega)$  
we define  the normal derivative of the $H^1$- solution $u$ of \eqref{eq:cauchy1} as ${\cal D_\nu}(u_0,u_1)$ and use 
the standard notation $\partial_\nu u$ instead of ${\cal D_\nu}(u_0,u_1)$ in such a way that we get \eqref{eq:hidden-alpha0}.
\end{Proof}
\begin{remark}
Theorem \ref{th:hidalpha} does not follow from the strong trace theorems of the Sobolev spaces. For this reason it can be called a hidden regularity result. The corresponding inequality \eqref{eq:hidden-alpha0}
is often called a direct inequality.
\end{remark}

%

\section{Conclusions and further questions}
In this paper we prove existence and hidden regularity for weak solutions of the fractional boundary value problem
\begin{equation*}
\begin{cases}
\displaystyle
\partial_t^{\alpha}u+\Delta^2 u=0
 \hskip0.9cm \mbox{in}\ (0,T)\times\Omega,
\\
u=\Delta u=0    \hskip1.4cm \mbox{on}\ (0,T)\times\partial\Omega.
\end{cases}
\end{equation*}
To conclude our analysis, we direct reader's attention to the following open problems.
\begin{enumerate}
\item Investigate  the hidden regularity in the case of Dirichlet--Neumann boundary conditions   
\begin{equation*}
\begin{cases}
\displaystyle
\partial_t^{\alpha}u+\Delta^2 u=0
 \hskip0.9cm \mbox{in}\ (0,T)\times\Omega,
\\
u=\partial_\nu u=0    \hskip1.4cm \mbox{on}\ (0,T)\times\partial\Omega,
\end{cases}
\end{equation*}
with $\alpha\in(1,2)$.
This study needs a new setting of the spaces and also 
 a trick in order to apply the multiplier method without using the integration by parts.
\item
A generalization of our research is dealing with
polyharmonic operator $\Delta^{2m}$  of order $2m$ 
with   boundary conditions 
$$u=\Delta u=\ldots=\Delta^{2m-1}u=0.$$
The corresponding  problem is
\begin{equation*}
\begin{cases}
\partial_t^{\alpha}u+\Delta^{2m}  u=0
 \hskip3.6cm \mbox{in}\ (0,T)\times\Omega,
\\
u=\Delta u=\ldots=\Delta^{2m-1}u=0$$
   \hskip1.4cm \mbox{on}\ (0,T)\times\partial\Omega.
\end{cases}
\end{equation*}
 To carry out this study, 
first one has to understand in which spaces the solution $u$ lives, and then to analyze how
  the methods used in this paper may be applied.

%
%
%
%

\end{enumerate}

\end{document}